\documentstyle[11pt]{article}

\newtheorem{definition}{Definition}
\newtheorem{theorem}{Theorem}
\newtheorem{lemma}{Lemma}
\newtheorem{proposition}{Proposition}
\newfont{\msbm}{msbm10}
\def\01{\mathop{[0,1]}\nolimits}
\def\Arc{\mathop{\rm Arc}\nolimits}
\def\E{\mathop{\mbox{\msbm E}}\nolimits}
\def\Geo{\mathop{\rm Geo}\nolimits}
\def\Lin{\mathop{\rm Lin}\nolimits}
\def\R{\mathop{\mbox{\msbm R}}\nolimits}
\def\V{\mathop{\mbox{\msbm V}}\nolimits}
\def\codom{\mathop{\rm codom}\nolimits}
\def\dom{\mathop{\rm dom}\nolimits}
\def\id{\mathop{\rm id}\nolimits}
\newcounter{example}
\newcommand{\example}[1]{\medskip\noindent\refstepcounter{example}%
{\bf Example \arabic{example} (#1). }}

\begin{document}

\title{The functional definition of geodesics}
\author{Lubom\'\i r Klapka}
\date{}

\maketitle

\begin{quote}\small
{\bf Abstract. } In this paper a functional definition of geodesics is
introduced
which allows to generalize the notion of a geodesic from smooth to
topological manifolds. It is shown that in the $C^\infty$-case the new
definition coincides with the classical definition of geodesics of a linear
connection. If $C^\infty$-smoothness is not required, it is shown by an
example that the new definition includes geodesics of non-linear, homogeneous
connections. Moreover, an example of generalized geodesics which does not
arise from any connection is presented here.
\end{quote}

\section{Introduction}

The concept of a geodesic is usually introduced by means of a linear
connection which has been developed from the Levi--Civita parallelism in
Riemannian geometry. The general geometric theory of linear and non-linear
connections, as applied in this paper, can be found in~\cite{Matsumoto}. In
this approach, a geodesic is always a solution of a second order differential
equation. In Sec.~\ref{1} we are recalling the well-known definition of a
linear connection based on the chart expression of this differential
equation, and some relations between geodesics and geodesic arcs.

Sec.~\ref{10} is devoted to open convex sets in a smooth manifold. Using the
Whitehead lemma~\cite{Whitehead}, we prove the existence of a covering of a
smooth manifold by open sets whose intersections are convex. We also derive
here functional equations~(\ref{16}), (\ref{17}) for geodesic arcs on a
convex manifold.

In Sec.~\ref{18}, these functional equations are generalized to manifolds
which are not necessarily convex and a new definition of geodesics is
introduced. Since the equations do not include derivatives, the new
definition makes sense on any topological manifold. We also give here two
important properties of the solutions of the functional equations.

In Sec.~\ref{24}, we prove that in the $C^\infty$-case, the functional
definition of a geodesic coincides with the usual one. The proof in one
direction uses the covering constructed in Sec.~\ref{10}. The proof of the
converse is based on the author's ideas contained in~\cite{Klapka}.

In Sec.~\ref{28}, we discuss some solutions of~(\ref{16}) and~(\ref{17})
which are obtained without the smoothness assumption. Since the general
solution is not known, a few particular examples are considered: a linear
solution, a continuous solution, a solution equivalent to a non-linear
connection, and a solution which does not correspond with any connection.

By a manifold we mean in this paper a finite-dimensional manifold of class
$C^k$, where $k\in\{0,1,2,\dots,\infty\}$; topologically we assume a
manifold to be second countable and Hausdorff. Manifolds with a boundary are
also used. For more information, the reader is referred to~\cite{Kobayashi}
and~\cite{Lang}.

Throughout this paper we will use the notation ${\partial f\over \partial
x^i}$ to mean $(f\circ x^{-1})_i\circ x$, where $f$ is a $C^1$-mapping,
$\codom f$ is a real finite-dimensional vector space, $x^i$ is the $i$th
coordinate function of a coordinate system $x$ on $\dom f$, and $(f\circ
x^{-1})_i$ is the partial derivative of $f\circ x^{-1}$ with respect to the
$i$th variable.

\section{The usual definition of geodesics}
\label{1}

Let $M$ be a $C^2$-manifold. A {\it linear connection $\Gamma$\/} on $M$ is a
set of functions $\Gamma^i_{jk}:\dom x\to\R$, $i,j,k\in\{1,2,\dots,\dim M\}$,
satisfying
\begin{equation}
\label{2}
\frac{\partial^2\overline x^i}{\partial x^j\partial
x^k}+\overline\Gamma^i_{lm}\frac{\partial\overline x^l}{\partial
x^j}\frac{\partial\overline x^m}{\partial x^k}=\frac{\partial\overline
x^i}{\partial x^l}\Gamma^l_{jk}
\end{equation}
on $\dom x\cap\dom\overline x$, where $x$ and $\overline x$ runs over all
coordinate systems on $M$. The functions $\Gamma^i_{jk}:\dom x\to\R$ are
called {\it components\/} of $\Gamma$ with respect to $x$. The manifold
$M$ is called a {\it domain of $\Gamma$\/} and is denoted by $\dom\Gamma$.
A linear connection $\Gamma$ is said to be {\it smooth\/} if $\dom\Gamma$ is
$C^\infty$-manifold and all components of $\Gamma$ are $C^\infty$-functions.

\begin{definition}
\label{3}
{\rm Let $\Gamma$ be a linear connection, $I\subset\R$ an open interval. A
$C^2$-map $g:I\to\dom\Gamma$ is called a} geodesic of $\Gamma$ {\rm if for
any $\tau\in I$ there exist coordinate systems $t$ on $I$ and $x$ on
$\dom\Gamma$ such that $\tau\in\dom
t$, $t=\id_{\dom t}$, $g(\dom t)\subset\dom x$, and
\begin{equation}
\label{4}
\frac{\partial^2x^i\circ g}{\partial t^2}+\bigl(\Gamma^i_{jk}\circ g\bigr)
 \, \frac{\partial x^j\circ g}{\partial t}\,\frac{\partial x^k\circ g}
{\partial t}=0,
\end{equation}
where $\Gamma^i_{jk}$ are components of $\Gamma$ with respect to $x$.
We denote the set of all geodesics of $\Gamma$ by $\Geo\Gamma$.}
\end{definition}

A map $h:\01\to\dom\Gamma$ is called a {\it geodesic arc\/} of $\Gamma$ if
there exists a geodesic $g\in \Geo\Gamma$ such that $\01\subset\dom g$,
$h=g|_{\01}$. We denote the set of all geodesic arcs of $\Gamma$ by
$\Arc\Gamma$.

\begin{proposition}
\label{5}
Suppose $\Gamma$ is a linear connection on $M$, $I\subset\R$ is an
open interval, $g:I\to M$ is a $C^2$-map. Then the following three statements
are equivalent:
\begin{enumerate}
\item\label{6} $g$ is a geodesic of $\Gamma$.
\item\label{7} For every two real numbers $\alpha,\beta\in I$, the map
\begin{equation}
\label{8}
h:\01\ni\gamma\to g\left((1-\gamma)\alpha+\gamma\beta\right)\in M
\end{equation}
is a geodesic arc of $\Gamma$.
\item\label{9} For every real number $\tau\in I$ there exist real numbers
$\alpha,\beta\in I$ such that $\alpha<\tau<\beta$ and~$(\ref{8})$ is a
geodesic arc of $\Gamma$.
\end{enumerate}
\end{proposition}
{\bf Proof.} Suppose the statement~\ref{6} holds. From Definition~\ref{3} it
follows that $g\circ p|_{ p^{-1}(I)}\in\Geo\Gamma$, where
$p:\R\ni\tau\to((1-\tau)\alpha+\tau\beta)\in\R$. Since $h=g\circ
p|_{p^{-1}(I)}|_{\01}$, we obtain the statement~\ref{7}.

Suppose the statement~\ref{7} holds. Since the interval $I$ is open, there
must exist real numbers $\alpha,\beta$ such that the statement~\ref{9} holds.

Suppose the statement~\ref{9} holds; then $g$ satisfies differential
equation~(\ref{4}) on an open neighborhood $(\alpha,\beta)$ of $\tau$.
Since $\tau\in I$ is arbitrary, $g$ satisfies this equation on $I$. Now, from
Definition~\ref{3} it follows that the statement~\ref{6} holds.

This completes the proof.

\section{Convex sets}
\label{10}

Let $\Gamma$ be a linear connection. An open set $U\subset\dom\Gamma$ is {\it
convex\/} with respect to $\Arc\Gamma$ if for any two points $a,b\in U$
there exists a unique geodesic arc $h\in\Arc\Gamma$ such that $h(0)=a$,
$h(1)=b$, $h(\01)\subset U$. Let us denote this geodesic arc by $h_{ab}$. An
open set $U\subset\dom\Gamma$ is {\it smoothly convex\/} with respect to
$\Arc\Gamma$ if $U$ is convex with respect to $\Arc\Gamma$ and the map
\begin{equation}
\label{11}
f: U\times U\times\01\ni(a,b,\tau)\to h_{ab}(\tau)\in U
\end{equation}
is smooth. It is clear that any open convex subset of a smoothly convex set
is smoothly convex. The empty set is, by definition, convex and smoothly
convex.
The existence of smoothly convex neighborhoods was proved by
J. H. C. Whitehead~\cite{Whitehead} in 1931. The proof of following lemma can
be also found in~\cite{Helgason} and~\cite{Kobayashi}.

\begin{lemma}[Whitehead]
\label{12}
Suppose $\Gamma$ is a smooth linear connection on a manifold $M$, $W\subset
M$ is an open set, $a\in W$ is a point; then there exists an open set
$V\subset M$ such that $a\in V\subset W$, and the set $V$ is smoothly convex
with respect to $\Arc\Gamma$.
\end{lemma}

The following proposition is needed for the sequel. Its proof is trivial.

\begin{proposition}
\label{13}
Suppose $\Gamma$ is a smooth linear connection, an open set
$U\subset\dom\Gamma$ is convex with respect to $\Arc\Gamma$,
$h,h'\in\Arc\Gamma$ are geodesic arcs such that $h(\01)\subset U$,
$h'(\01)\subset U$, $h(0)=h'(0)$, $h(1)=h'(1)$. Then $h=h'$.
\end{proposition}

We note that convex sets in a Euclidean space have different properties than
convex sets in a manifold endowed with a linear connection. In particular,
the intersection of two convex sets of a Euclidean space is always convex,
while the intersection of two convex sets in the circle is not necessarily
convex, and even may be non-connected. Whitehead's lemma guarantees existence
of a covering by open, convex sets; nothing is claimed, however, about their
intersections. In what follows, we need the following result.

\begin{lemma}
\label{14}
If $\Gamma$ is a smooth linear connection on manifold $M$, then there exists
an open covering $\{U_a\}_{a\in M}$ of $M$, where $U_a$'s are neighborhoods
of $a$'s, such that the set $U_a\cap U_b$ is smoothly convex with respect to
$\Arc\Gamma$ for every two points $a,b\in M$.
\end{lemma}
{\bf Proof.} Consider a smooth Riemannian metric on $M$. Let $a\in M$ be a
point, $x$ a coordinate system on $M$ such that $a\in\dom x$. Let us consider
the function $u:\dom x\ni b\to\triangle^2(a,b)\in\R$, where $\triangle(a,b)$
is the distance of $b$ from $a$. At the point $a$, the first partial
derivatives of $u$ vanish, and the second partial derivatives are components
of the metric tensor. This has been proved by J. L. Synge in~\cite{Synge}.
Therefore, the second covariant derivatives $\partial^2u/\partial x^k\partial
x^l-\Gamma^m_{kl}\partial u/\partial x^m$ at $a$ are coefficients in a
positive definite quadratic form. From the continuity it follows that there
exists an open set $W$ such that $a\in W\subset\dom x$ and the second
covariant derivatives of $u$ are coefficients in a positive definite
quadratic form at each point $b\in W$. We restrict $x$ to $W$.

By Lemma~\ref{12}, there exists an open set $V$ such that $a\in V\subset\dom
x$ and $V$ is smoothly convex with respect to $\Arc\Gamma$. Further, there
exists an open ball $B_a\subset V$ of center $a\in B_a$ and radius $\rho_a$.
Let us consider two points $c,d\in B_a$. Since $V$ is convex, there is a
unique geodesic arc $h\in\Arc\Gamma$ such that $h(0)=c$, $h(1)=d$, and
$h(\01)\subset V$. From~(\ref{4}), we get
$$
\frac{\partial^2u\circ h}{\partial t^2}=\left(\Bigl(\frac{\partial^2u}
{\partial x^k\partial x^l}-\Gamma^m_{kl}\frac{\partial u}{ \partial x^m}
\Bigr)\circ h\right)\frac{\partial x^k\circ h}{\partial t}\frac{\partial
x^l\circ h}{\partial t}.
$$
Since the quadratic form in the first derivatives is positive definite, we
have $\partial^2u\circ h/\partial t^2\ge0$, whence the function $u\circ h$ is
convex. In this case, for each $\tau\in\01$, we have
$$
u\circ h(\tau)\le(1-\tau)u\circ h(0)+\tau u\circ h(1)<(1-\tau)
\rho^2_a+\tau\rho^2_a=\rho^2_a.
$$
This means that $h(\01)\subset B_a$. Thus the open set $B_a\subset M$ is
smoothly convex with respect to $\Arc\Gamma$. For the same reason, any open
ball of center $a$ and radius smaller than $\rho_a$ is smoothly convex with
respect to $\Arc\Gamma$.

Thus, there exists an open covering $\{B_a\}_{a\in M}$ of $M$, where
any $B_a$ is a smoothly convex open ball of center $a$ and radius $\rho_a$.
Moreover, there exists an open covering $\{U_a\}_{a\in M}$, where any
$U_a$ is the smoothly convex open ball of center $a$ and radius
$\frac{1}{3}\rho_a$.

Now, consider two points $a,b\in M$. For the sake of being definite assume
that $\rho_a\le\rho_b$. If $U_a\cap U_b=\emptyset$, then the open set
$U_a\cap U_b$ is smoothly convex with respect to $\Arc\Gamma$. Conversely,
suppose that $U_a\cap U_b\not=\emptyset$. Consider two points $c,d\in U_a\cap
U_b$. Since $U_a$, $U_b$ are convex sets, there exist geodesic arcs
$h,h'\in\Arc\Gamma$ such that $h(0)=h'(0)=c$, $h(1)=h'(1)=d$, $h(\01)\subset
U_a$, and $h'(\01)\subset U_b$. From $\frac{2}{3}\rho_a+\frac{1}{3}\rho_b\le
\rho_b$ it follows that $U_a\cup U_b\subset B_b$. Since $B_b$ is a convex
set, by Proposition~\ref{13}, $h=h'$. Hence, we get $h(\01)=h'(\01)\subset
U_a\cap U_b$. From here, the open set $U_a\cap U_b$ is, by definition,
smoothly convex with respect to $\Arc\Gamma$.

This completes the proof of Lemma~\ref{14}.

In the following proposition, we consider functional properties of the
mapping~(\ref{11}).

\begin{proposition}
\label{15}
Suppose $\Gamma$ is a linear connection on a manifold $M$, $U\subset M$ is an
open set, convex with respect to $\Arc\Gamma$, $f$ is the
mapping~$(\ref{11})$, $a,b\in U$ are any points, and
$\alpha,\beta,\gamma\in\01$ are real numbers. Then
\begin{equation}
\label{16}
f(a,b,0)=a,\quad f(a,b,1)=b,
\end{equation}
\begin{equation}
\label{17}
f(a,b,(1-\gamma)\alpha+\gamma\beta)=f(f(a,b,\alpha),f(a,b,\beta),\gamma).
\end{equation}
\end{proposition}
{\bf Proof.} By construction, $f$ satisfies~(\ref{16}), $\01\ni\tau\to
f(a,b,\tau)\in M$ is a geodesic arc of $\Gamma$. From Proposition~\ref{5},
the mappings $h:\01\ni\tau\to f(a,b,(1-\tau)\alpha+\tau\beta)\in M$,
$h':\01\ni\tau\to f(f(a,b,\alpha),f(a,b,\beta),\tau)\in M$ are geodesic arcs
of $\Gamma$. Combining~(\ref{16}) and Proposition~\ref{13}, we
obtain~(\ref{17}).

This completes the proof.

\section{The new definition of geodesics}
\label{18}

Let $M$ be a finite-dimensional topological manifold. In this section, we
generalize equations~(\ref{16}) and~(\ref{17}) from convex smooth manifolds
to non-convex topological manifolds. For this purpose, we replace the smooth
mapping~(\ref{11}) by a continuous mapping
\begin{equation}
\label{19}
f: D\times\01\to M,
\end{equation}
where $D\subset M\times M$ is an open set containing the diagonal, and
(\ref{16}), (\ref{17}) hold for each $(a,b)\in D$,
$\alpha,\beta,\gamma\in\01$.

By Acz\'el definitions~\cite{Aczel}, the formulas~(\ref{16}), (\ref{17}) are
called {\it functional equations\/}, and a mapping~(\ref{19}) is called a
{\it solution of functional equations\/}. A solution of functional equation
is said to be {\it smooth\/} if (\ref{19}) is $C^\infty$.

\begin{definition}
\label{20}
{\rm A mapping $g:I\to M$, where $I\subset\R$ is an open interval, is called
a} geodesic, {\rm associated with a solution~(\ref{19}) of the functional
equations~(\ref{16}), (\ref{17}) if and only if to each $\tau\in I$ there
exist $\alpha,\beta\in I$, $\alpha<\tau<\beta$, such that for each
$\gamma\in\01$,
\begin{equation}
\label{21}
g\left((1-\gamma)\alpha+\gamma\beta\right)=f\left(g(\alpha),g(\beta),\gamma
\right).
\end{equation}
The set of all geodesics associated with $f$ is denoted by $\Geo f$.}
\end{definition}

An easy consequence of the definition is that geodesics associated with a
solution of the functional equations are always continuous.

Before going on, we wish to mention two special cases. Let $f$ be a solution
of~(\ref{16}) and~(\ref{17}), let $(a,b,\gamma)\in\dom f$. Taking $\alpha
=\beta =0$ in~(\ref{17}) we get by~(\ref{16})
\begin{equation}
\label{22}
f(a,a,\gamma)=a.
\end{equation}
Taking $\alpha=1$, $\beta=0$, we obtain
\begin{equation}
\label{23}
f(a,b,1-\gamma)=f(b,a,\gamma).
\end{equation}

\section{The comparison of the definitions in $C^\infty$-case}
\label{24}

In this section, we wish to prove that in the $C^\infty$-case, the new
definition of geodesics coincides with the usual one.

\begin{theorem}
To every smooth linear connection $\Gamma$ there exists a smooth solution $f$
of the functional equations~$(\ref{16})$ and $(\ref{17})$ such that
$\Geo f=\Geo\Gamma$.
\end{theorem}
{\bf Proof.} 1. First, we construct a solution $f$. For this purpose we
consider a covering $\{U_a\}_{a\in M}$ of $M=\dom\Gamma$ with the properties
described by Lemma~\ref{14}. To each $U_a$ we assign a smooth solution $f_a:
U_a\times U_a\times\01\to U_a$ of~(\ref{16}), (\ref{17}) (see
Proposition~\ref{15}).

Let $a,b\in M$, $c,d\in U_a\cap U_b$. Since $U_a\cap U_b$ is convex, there
exists just one geodesic arc $h\in\Arc\Gamma$ such that $h(0)=c$,
$h(1)=d$, $h(\01)\subset U_a\cap U_b$. Since $U_a$, $U_b$ are convex, we get
for all $\gamma\in\01$ the equality
$h(\gamma)=f_a(c,d,\gamma)=f_b(c,d,\gamma)$. If $c,d\in U_a\cap U_b$ are
arbitrary, then $f_a|_{\dom f_a\cap\dom f_b}=f_b|_{\dom f_a\cap\dom f_b}.$
Moreover, since $a,b\in M$ are arbitrary, there should exist a mapping
$f:D\times\01\to M$, where
$$
D=\bigcup_{a\in M} U_a\times U_a\subset M\times M,
$$
such that for each $a\in M$, $f|_{\dom f_a}=f_a$. The set $D$ is open and
contains the diagonal of $M\times M$. Since $f_a$ are smooth solutions of
(\ref{16}), (\ref{17}), $f$ is also a smooth solution of~(\ref{16}),
(\ref{17}).

2. Now we prove that $\Geo f=\Geo\Gamma$. Let $g\in\Geo f$. Then by
Definition~\ref{20}, to each $\tau\in\dom g$ there exist $\alpha,\beta\in\dom
g$ such that relations $\alpha<\tau<\beta$ and~(\ref{21}) hold. The
equality~(\ref{21}) can be written as $h=h'$, where $h$ is defined
by~(\ref{8}), and $h'$ by
\begin{equation}
\label{25}
h':\01\ni\gamma\to f\left(g(\alpha),g(\beta),\gamma\right)\in M.
\end{equation}
The construction of~$f$ implies $h'\in\Arc\Gamma$, and $\tau\in\dom g$ is
arbitrary. Thus, Proposition~\ref{5} implies $g\in\Geo\Gamma$.

Conversely, assume that $g\in\Geo\Gamma$, $\tau\in\dom g$, $a=g(\tau)$. Since
$g$ is continuous, there exist $\alpha,\beta\in\dom g$ for which
$\alpha<\tau<\beta$, $g([\alpha,\beta])\subset U_a$. By proposition~\ref{5},
$h\in\Arc\Gamma$, and by construction of~$f$, $h'\in\Arc\Gamma$. Using
convexity of $U_a$, (\ref{16}) and Proposition~\ref{13}, we get the
relation~(\ref{21}). Since $\tau\in\dom g$ is arbitrary, Definition~\ref{20}
implies $g\in\Geo f$.

This completes the proof.

\begin{theorem}
To every smooth solution $f$ of the functional equations~$(\ref{16})$ and
$(\ref{17})$ there exists a smooth linear connection $\Gamma$ such that $\Geo
f=\Geo\Gamma$.
\end{theorem}
{\bf Proof.} 1. Let us construct $\Gamma$. Let $f: D\times\01\to M$ be a
smooth solution of~(\ref{16}) and~(\ref{17}), $a\in M$ an arbitrary point,
$x$ a coordinate system at $a$ such that $\dom x\times\dom x\subset D$. Let
$n=\dim M$. Denote
\begin{equation}
\label{26}
\Gamma ^i_{jk}(a)=-4\,\frac{\partial ^2 x^i\circ f}{\partial y^j\partial
z^k}\,(a,a,{\textstyle\frac{1}{2}}),
\end{equation}
where $y^1,y^2,\dots,y^n,z^1,z^2,\dots,z^n,t$ are the coordinate functions of
the coordinate system $x\times x\times\id_{\01}$ on $D\times\01$ (compare
with~\cite{Klapka}). By~(\ref{22}) and~(\ref{23}),
\begin{equation}
\label{27}
\frac{\partial x^i\circ f}{\partial
y^j}\,(a,a,{\textstyle\frac{1}{2}})=\frac{\partial x^i\circ f}{\partial
z^j}\,(a,a,{\textstyle\frac{1}{2}})={\textstyle\frac{1}{2}}\delta^i_j,
\end{equation}
where $\delta^i_j$ is the Kronecker symbol. Since $a\in M$ is arbitrary,
(\ref{22}), (\ref{26}), and (\ref{27}) imply~(\ref{2}). Thus, (\ref{26})
defines a linear connection $\Gamma$ on $M$. Since the components of $\Gamma$
are restrictions of the second derivatives of smooth functions to a smooth
submanifold, $\Gamma$ must be smooth.

2. Now we prove that $\Geo f=\Geo\Gamma$. Let $g\in\Geo f$, $\tau\in\dom g$.
Then by Definition~\ref{20} there exist $\alpha_0,\beta_0\in\dom g$ such that
for $\alpha=\alpha_0$, $\beta=\beta_0$, and $\gamma\in\01$, relations
$\alpha<\tau<\beta$ and ~(\ref{21}) hold. Then it follows for
$\alpha,\beta\in[\alpha_0,\beta_0]$, $\gamma\in\01$ that
$$
g\left((1-\gamma)\alpha+\gamma\beta\right)=f\left(g(\alpha_0),g(\beta_0),
\frac{(1-\gamma)(\alpha-\alpha_0)+\gamma(\beta-\alpha_0)}{\beta_0-\alpha_0}
\right).
$$
But then by~(\ref{17}), the relation~(\ref{21}) holds for any
$\alpha,\beta\in[\alpha_0,\beta_0]$, $\gamma\in\01$. Consider coordinate
systems $t$ on $\dom g$ and $x$ on $M$ such that $\tau\in\dom t$,
$t=\id_{\dom t}$, $g(\dom t)\subset\dom x$, and $\dom x\times\dom x\subset
D$. Further, consider the chart $x\times x\times\id_{\01}$ on $D\times\01$.
Differentiating the coordinate expression of~(\ref{21}) with respect to
$\alpha$ and $\beta$ and setting $\alpha=\beta=\tau$, $\gamma=\frac{1}{2}$,
we find, using~(\ref{26}), that $g$ satisfies~(\ref{4}) at the point $\tau$.
According to~(\ref{2}), the equation~(\ref{4}) is satisfied also for $\dom
x\times\dom x\not\subset D$. Since $\tau\in\dom g$ is arbitrary, and the
charts $t$, $x$, such that $\tau\in\dom t$, $t=\id_{\dom t}$, $g(\dom
t)\subset\dom x$, are also arbitrary, we get, using Definition~\ref{3}, that
$g$ is a geodesic of the connection $\Gamma$, i.e., $g\in\Geo\Gamma$.

To prove the converse, choose $g\in\Geo\Gamma$, $\tau\in\dom g$. By
Lemma~\ref{12}, there exists a neighborhood $U\subset M$ of the point
$g(\tau)$, convex with respect to $\Arc\Gamma$. Then by~(\ref{22}), the
continuity of $f$, and the compactness of the interval $\01$, imply existence
of an open set $W\subset M$ such that $g(\tau)\in W$, $W\times
W\times\01\subset\dom f$, and $f(W\times W\times\01)\subset U$. Consider the
mappings $h$ and $h'$ defined by~(\ref{8}) and~(\ref{25}). Since $g$ is
continuous, the numbers $\alpha,\beta\in\dom g$ can be chosen in such a way
that $\alpha<\tau<\beta$ and $h(\01)\subset U\cap W$. Then indeed
$h'(\01)\subset U$. By Proposition~\ref{5}, $h\in\Arc\Gamma$, and
by~(\ref{17}) and Definition~\ref{20}, $h'|_{(0,1)}\in\Geo f$. But we already
know that $\Geo f\subset\Geo\Gamma$, hence $h'|_{(0,1)}$ is a solution
of~(\ref{4}). Since $f$ is smooth, the domain of the definition of this
solution can be extended to a neighborhood of $\01$, hence $h'\in\Arc\Gamma$.
The convexity of $U$ implies, by~(\ref{16}) and Proposition~\ref{13},
that~(\ref{21}) holds. Now since $\tau\in\dom g$ is arbitrary, we have, using
Definition~\ref{20}, $g\in\Geo f$, as required.

This completes the proof.

\section{Examples}
\label{28}

The following examples are obtained by a reparametrization of the geodesic
arcs corresponding with the linear solution of functional equations.

\example{linear solution}
\label{29}
Let $\V$ be a real, finite-dimensional vector space. We look for a continuous
solution $f:\V\times \V\times\01\to \V$ of~(\ref{16}) and (\ref{17}) such
that for each $\tau\in\01$, the mapping $\V\times \V\ni (a,b)\to
f(a,b,\tau)\in \V$ is linear.

By linearity, there exists a mapping $q:\01\to\Lin(\V,\V)$, assigning to
every $\gamma\in\01$ a linear mapping $q(\gamma):\V\ni b\to f(0,b,\gamma)\in
\V$. By~(\ref{16}),
\begin{equation}
\label{30}
q(0)=0.
\end{equation}
Since for every $a,b\in\V$, $\gamma\in\01$ we get with help of~(\ref{23})
$f(a,b,\gamma)=q(1-\gamma)(a)+q(\gamma)(b)$, it follows from~(\ref{22}) that
\begin{equation}
\label{31}
q({\textstyle\frac{1}{2}})={\textstyle\frac{1}{2}}\id_{\V}.
\end{equation}
Substituting into~(\ref{17}) $a=0$, $\gamma=\frac{1}{2}$, we get,
using~(\ref{23}) and~(\ref{31}) the Jensen's functional equation
$$
q\left(\frac{\alpha+\beta}{2}\right)=\frac{q(\alpha)+q(\beta)}{2}.
$$
Under conditions~(\ref{30}) and~(\ref{31}), this equation has a unique
continuous solution
$q(\gamma)=\gamma\id_{\V}$ (see e.g.~\cite{Aczel}). Thus,
$$
f(a,b,\gamma)=(1-\gamma)a+\gamma b,
$$
which satisfies~(\ref{16}) and~(\ref{17}).

It is easily seen that all geodesics associated with this solution are also
geodesics of the canonical linear connection on the manifold $\V$.

\example{continuous solution}
\label{32}
Let $\E$ be a real, finite-dimensional vector space,
$\E\times\E\ni(a,b)\to(ab)\in\R$ a scalar product, $\E\ni
a\to|a|\in[0,\infty)$ the corresponding norm, and let $v:\R\to\R$ be a
homeomorphism such that for each $\alpha\in\R$, $v(-\alpha)=- v(\alpha)$. We
shall show that the formula
\begin{equation}
\label{33}
f(a,b,\gamma)=\cases{ v_{ab}^{-1}\bigl((1-\gamma) v_{ab}(a)+\gamma v_{ab}(b)
\bigr)&for $a\not=b$\cr a&for $a=b$, \cr}
\end{equation}
where $v_{ab}:\E\to\E$ is defined by $v_{ab}(c)=c+(v((ce))-(ce))e$,
$e=(a-b)/|a-b|$, defines a continuous solution $f:\E\times\E\times\01\to\E$
of~(\ref{16}), (\ref{17}).

It is easy to verify that if $a\not=b$, then $v_{ab}$ is a bijection, and
$v_{ab}^{- 1}(c)=c+( v^{-1}((ce))-(ce))e$. This shows the existence
of~(\ref{33}). If $a=b$, or $a\not=b$ and $\alpha=\beta$, $f$ obviously
satisfies~(\ref{16}) and~(\ref{17}). Consider the remaining case $a\not=b$,
$\alpha\not=\beta$. By a straightforward computation,
$$
\frac{f(a,b,\alpha)-f(a,b,\beta)}{|f(a,b,\alpha)-f(a,b,\beta)|}=\frac{a-b}
{|a-b|}\mbox{sign}(\beta-\alpha).
$$
Therefore, $v_{f(a,b,\alpha)f(a,b,\beta)}=v_{ab}$ and $f$
satisfies~(\ref{16}), (\ref{17}) as required.

Now, we prove continuity of $f$. For all $a,b\in\E$, $\gamma\in\01$, it
follows that $|a-f(a,b,\gamma)|+|f(a,b,\gamma)-b|=|a-b|$. By the triangle
inequality, the point $f(a,b,\gamma)$ belongs to the segment joining $a$ and
$b$. Thus, for every open ball $B\subset\E$, $f(B\times B\times\01)=B$. This
proves the convexity of all open balls with respect to $\Arc f$, as well as
the continuity of $f$ for $a=b$. If $a\not=b$, the continuity of $f$ is
obvious.

\example{a solution equivalent to a non-linear connection}
Consider Example~\ref{32} and suppose moreover that $v$ is a diffeomorphism
of class $C^\infty$. Since for $a=b$ some of the second partial derivatives
of~(\ref{33}) do not exist, $f$ is of class $C^1$. Let $r:\R\to\R$ be a
mapping defined by
$$
\frac{\partial^2 v}{\partial t^2}=r\frac{\partial v}
{\partial t},
$$
where $t=\id_{\R}$, and $s:\E\times\E\to\E$ is defined by
$$
s(a,b)=\cases{r\left((ab)/|b|\right)|b|\,b&for $b\not=0$\cr 0&for $b=0$.\cr}
$$
Then all geodesics $g\in\Geo f$ are solutions of the differential equation
$$
\frac{\partial^2g}{\partial t^2}+s\left(g,\frac{\partial g}{\partial
t}\right)=0.
$$
Since for each $\lambda\in\R$, $a,b\in\E$, we have $s(a,\lambda
b)=\lambda^2s(a,b)$, these geodesics can be interpreted as geodesics of a
non-linear, homogeneous connection (see e.g.~\cite{Matsumoto}). In some
special cases, e.g. for $v$ defined by
$$
v:\R\ni\alpha\to\int\limits^\alpha_0\exp\left(\frac{\beta^2}{2}\right)d
\beta\in\R,
$$
this connection is linear and smooth. If $v=\id_{\R}$ we get the situation
described by Example~\ref{29}.

\example{a solution which does not represent a connection}
Let $w:\R\to\R$ be the Weierstrass's function; $w$ is continuous and bounded,
but it is not differentiable at any point (see~\cite{Hardy}). Since $w$ is
bounded, there exists $\kappa\in(0,\infty)$ such that for each $\alpha\in\R$,
$-\kappa<w(\alpha)<\kappa$. Then by the continuity of $w$, the function
\begin{equation}
\label{34}
\R\ni\alpha\to\int\limits^\alpha_0\left(w(\beta)+\kappa\right)d\beta\in\R,
\end{equation}
is differentiable on $\R$, and its derivative is positive. Since the
indefinite integral of the Weierstrass's function is a bounded function on
$\R$, it follows that (\ref{34}) is a surjection. Hence, (\ref{34}) is a
diffeomorphism of class $C^1$.

Let $v^{-1}$ in Example~\ref{32} be the mapping~(\ref{34}). Then
by~(\ref{21}) and~(\ref{33}), the geodesic $g\in\Geo f$ is not twice
differentiable at any point of its domain. Therefore, $g$ cannot be a
solution of a second order differential equation. Consequently, there is
neither a linear nor non-linear connection whose geodesic is $g$.

\section*{Acknowledgments}

The author is grateful to professor Demeter Krupka for his constant attention
to this work. He also thanks Jaroslav \v Stef\'anek and Michal Marvan for
useful discussions. This work has been supported from the Project VS 96003
``Global Analysis'' by the Czech Ministry of Education, Youth, and Sports.

\footnotesize
\noindent Lubom\'{\i}r Klapka \\
Department of Mathematics and Computer Science \\
            Silesian University at Opava \\
            Bezru\v{c}ovo n\'am.~13 \\
            746~01 Opava 1, Czech Republic

\end{document}